\documentclass[ECP]{ejpecp} 

\AUTHORS{%
  Xiequan Fan\footnote{Center for Applied Mathematics,
Tianjin University, Tianjin 300072,  China  \ \ \ \ \ \newline  \mbox{\ \ \ \ \ }   E-mail: fanxiequan@hotmail.com}
     }

\SHORTTITLE{On McDiarmid's  inequality}

\TITLE{On McDiarmid's  inequality for  Hamming distance} 


\AMSSUBJ{60E15} 

\SUBMITTED{September  ., 2016} 
\ACCEPTED{December ., 2016} 

\VOLUME{0}
\YEAR{2016}
\PAPERNUM{0}
\DOI{vVOL-PID}

\ABSTRACT{We improve the rate function of  McDiarmid's inequality for Hamming distance.
In particular, applying our result to the separately Lipschitz functions of independent random variables,
we also refine  the convergence rate function of McDiarmid's inequality around a median.
Moreover, a non-uniform bound for the distance between  the medians and the mean is also given.
We also give some extensions of  McDiarmid's inequalities to the case of nonnegative functionals of  dependent  random variables.}

\begin{document}



\section{Introduction}
Let $\alpha=(\alpha_1,...,\alpha_n) \in [0, \infty)^n$ be an $n$-vector of non-negative real numbers. Recall that the
$L_2$ norm is given by
\[
||\alpha||=(\sum_{i=1}^n \alpha_i^2)^{1/2},
\]
and we call $\alpha$ a unit vector if $||\alpha||=1.$ For points $x=(x_1,...,x_n)$ and $y=(y_1,...,y_n)$ in $\mathbf{E}^n$,
 the $\alpha$-Hamming distance $d_\alpha(x, y)$ is defined as
 \[
 d_\alpha(x, y) = \sum_{x_i\neq y_i} \alpha_i.
 \]
In particular, when $\alpha_i=1/\sqrt{n},$ the $\alpha$-Hamming distance is just the unit Hamming distance. For
a subset $A \subset \mathbf{E}^n,$ we define
 \[
 d_\alpha(x, A) = \inf_{y \in A} d_\alpha(x, y).
 \]
All over the paper, we assume that $\mathbf{E}$ is a separable
spaces such that $d_\alpha(x, A)$ is measurable.

McDiarmid proved the following concentration inequality for $\alpha$-Hamming distance   (see Theorem 3.6 of \cite{M98}).
Let $X$ be an $\mathbf{E}^n$-valued vector of independent random variables.
Assume that $||\alpha||=1.$ Then, for any set $A \subset \mathbf{E}^n$
and any $t>0,$
\begin{eqnarray}\label{McD01}
\mathbb{P}(d_\alpha(X, A)  \geq t) \,  \mathbb{P}( X \in A)  \leq e^{ - t^2/2  }.
\end{eqnarray}
Such result is very useful when one wants to evaluate the concentration around a median.
To illustrate  its application, consider  the following example. All over the paper, let $f $ be a function defined on $\mathbf{E}^n$.
Assume  that
the function $f$ satisfies
\begin{eqnarray}\label{McD02}
 |f(x)-f(x')| \leq d_\alpha(x, x') .
\end{eqnarray}
Let $m$ be a median of $f(X),$ that is $\mathbb{P}(f(X) \geq m  ) \geq \frac12 $  and $\mathbb{P}(f(X) \leq m  ) \geq \frac12 $. By taking $ A=\{ X\in \mathbf{E}^n: f(X)\leq m\},$   McDiarmid's inequality (\ref{McD01}) implies the following concentration inequalities around a median: for any  $t>0,$
 \begin{eqnarray}\label{ghlm1}
\mathbb{P}(f(X) -m  \geq t)   \leq 2\, e^{ -\frac12 t^2  }
\end{eqnarray}
and
 \begin{eqnarray}\label{ghlm2}
\mathbb{P}(f(X) -m  \leq - t)   \leq 2\, e^{ -\frac12  t^2   }.
\end{eqnarray}
Denote by $\mu $ the mean value of the random variable $f(X)$, that is  $\mu= \mathbb{E} [f(X)].$
By Lemma 4.6 of McDiarmid \cite{M98}, the  inequalities (\ref{ghlm1}) and (\ref{ghlm2})  together implies
a  uniform bound between  the medians  and the mean:
\begin{eqnarray}\label{fsg21}
  |m-\mu | \leq \sqrt{2 \pi}.
\end{eqnarray}

McDiarmid (see  Theorem 3.1 of \cite{M98}) also proved the following exponential moment inequality:
for any $\lambda\geq0,$
 \begin{eqnarray}\label{fsg250}
\mathbb{E}[e^{ \lambda(f(X) - \mu)}  ]   \leq \exp\Big\{ \frac{\lambda^2} 8 \Big \}.
\end{eqnarray}
The last inequality together with the classical Bernstein  inequality $$\mathbb{P}(f(X) - \mu  \geq t) \leq \inf_{\lambda \geq 0}e^{-\lambda t}\mathbb{E}[e^{ \lambda(f(X) - \mu)}  ]$$  implies   the following
 concentration inequality around the mean (instead of a median): for any  $t>0,$
 \begin{eqnarray}\label{fsg251}
\mathbb{P}(f(X) - \mu  \geq t)   \leq   e^{ - 2\, t^2   }.
\end{eqnarray}
 See also Rio \cite{R13} for a recent improvement of (\ref{fsg251}), such that $\mathbb{P}(f(X) - \mu  \geq t) =0$ for any
$t> \sum_{i=1}^n \alpha_i .$
By (\ref{fsg251}), we get  for any  $t>|m-\mu|,$
 \begin{eqnarray}\label{f26sf}
\mathbb{P}(f(X) - m \geq t)   \leq    e^{ - 2\, (t- |m-\mu|) ^2   }.
\end{eqnarray}
The last inequality shows that for large $t$'s, the constant $\frac12$ in the inequalities (\ref{ghlm1}) and (\ref{ghlm2}) can be improved.
However, a deficiency of (\ref{f26sf}) is that comparing to (\ref{ghlm1}) and (\ref{ghlm2}),  inequality (\ref{f26sf})   does not hold for   $t \in (0, |m-\mu|].$
Moreover, in applications we need to give a  bound for  $|m-\mu|$. Obviously, this  bound should be as small as possible.

The scope of the paper is to fill the deficiency of (\ref{f26sf}), and establish a bound sharper than (\ref{fsg21}).  To achieve this scope,  we give an improvement on (\ref{McD01}):
 for any set $A \subset \mathbf{E}^n$
and any $t >0,$
\begin{eqnarray}\label{fghl}
\mathbb{P}(d_\alpha(X, A)  \geq t) \,  \mathbb{P}( X \in A)  \leq e^{ -  t^2 }.
\end{eqnarray}
  It is obvious that the r.h.s.\ of   (\ref{fghl}) is much smaller than the one of (\ref{McD01})  due to the fact that the constant $\frac12$ has been
  improved to $1$.
Using (\ref{fghl}),  we obtain an improvement on (\ref{ghlm1}):  for any $t >0,$
 \begin{eqnarray}\label{ghds}
\mathbb{P}(f(X) -m  \geq t)   \leq 2\, e^{ -  t^2 }.
\end{eqnarray}
We also prove the following bound for  $|m-\mu|$:
\begin{eqnarray}\label{f10.8}
| \mu-m |  \leq  \big(2  \rho + \sqrt{\pi/2  } \, \big)e^{ -  2 \rho^2  }  \, <\, 1.5503,
\end{eqnarray}
where  $\rho= \mathbb{E}[d_\alpha(X, A)]$ with $A=\{Y: f(Y)\leq m \} $.
 Since $ 1.5503 < \sqrt{2 \pi}, $ inequality (\ref{f10.8}) is tighter than (\ref{fsg21}).
In particular, (\ref{f10.8})  implies a non-uniform bound between the  medians  and the mean, which shows that if $\rho \rightarrow \infty,$ then the distance between the  medians  and the mean tends to $0$ in an exponentially decaying rate.

This paper is organized as follows. Our main results are stated and discussed in Section \ref{Sec2}.
In Section \ref{sec3}, we extend the McDiarmid  inequalities  to the case  of  nonnegative functionals of \emph{dependent} random variables.
We prove our theorems   in Section \ref{Sec3}.

\section{On McDiarmid's  inequalities}\label{Sec2}
In the following theorem, we give a stronger version of McDiarmid's inequality (\ref{McD01}).
\begin{theorem}\label{th1}  Let $X$ be an $\mathbf{E}^n$-valued vector of independent random variables.
Assume that $||\alpha||=1.$   Then for any set $A \subset \mathbf{E}^n$
and any $t>0,$
\begin{eqnarray}\label{sf01}
 \mathbb{P}(d_\alpha(X, A)  \geq t) \,  \mathbb{P}( X \in A) \leq  e^{ - h(t)   },
\end{eqnarray}
where
\begin{eqnarray} \label{sf02}
 h(t)=   \left\{ \begin{array}{ll}
2 (  \mathbb{E}[d_\alpha(X, A)])^2, & \textrm{\ \ \ \ \ if $0< t< \mathbb{E}[d_\alpha(X, A)]$,}\\
\\
t^2+ (t-2 \mathbb{E}[d_\alpha(X, A)])^2, & \textrm{\ \ \ \ \  if $t\geq \mathbb{E}[d_\alpha(X, A)]$.}\\
\end{array} \right.
\end{eqnarray}
In particular, it implies that for any $t>0,$
\begin{eqnarray}\label{ine10}
\mathbb{P}(d_\alpha(X, A)  \geq t) \,  \mathbb{P}( X \in A)
 \leq e^{ - t^2},
\end{eqnarray}
and that
\begin{eqnarray} \label{fdgs}
\mathbb{P}(X \notin A) \,  \mathbb{P}( X \in A)
 \leq e^{ -  2(  \mathbb{E}[d_\alpha(X, A)])^2}.
\end{eqnarray}
\end{theorem}

Inequality (\ref{sf01}) shows that the   constant $\frac12$ in (\ref{McD01}) can be improved  to $2$ for small $t$. Indeed,
when   $t \in (0, \,   \mathbb{E}[d_\alpha(X, A)]),$ it holds $ 2 ( \mathbb{E}[d_\alpha(X, A)])^2 \geq 2\,t^2. $
Thus the bound (\ref{sf01}) implies that for $t \in (0, \,    \mathbb{E}[d_\alpha(X, A)]),$
\begin{eqnarray}
\mathbb{P}(d_\alpha(X, A)  \geq t) \,  \mathbb{P}( X \in A)  \leq e^{ - 2 \, t^2  },
\end{eqnarray}
which concludes ours claim.

 When   $t\rightarrow \infty,$ it holds that $(t^2+  (t-2 \mathbb{E}[d_\alpha(X, A)])^2 )/(2\,t^2) \rightarrow 1. $
Thus the bound  (\ref{sf01}) also behaviors as $\exp\big\{- 2\, t^2  \big\}$ for large $t$'s.

 Inequality (\ref{ine10}) shows that the constant $\frac12$ in McDiarmid's inequality (\ref{McD01}) can be improved to $1$. The virtue of
   the bound (\ref{ine10}) is that it does not depend on $\mathbb{E}[d_\alpha(X, A)]$.

 It is well-known that  $\mathbb{P}(X \notin A) \,  \mathbb{P}( X \in A) \leq \frac14 \big(\mathbb{P}(X \notin A) +  \mathbb{P}( X \in A) \big)^2=\frac14.$
While, inequality (\ref{fdgs}) implies a stronger result: if $\mathbb{E}[d_\alpha(X, A)] \rightarrow \infty,$ then it holds
\begin{eqnarray*}
\mathbb{P}(X \notin A) \,  \mathbb{P}( X \in A)  \rightarrow 0.
\end{eqnarray*}
This  new feature   does not imply by  (\ref{McD01}).

Next, we apply Theorem \ref{th1}   to the study of concentration around a median. We have the following improvement on (\ref{ghlm1}) and (\ref{ghlm2}).
\begin{theorem}\label{th2}  Let $X$ be an $\mathbf{E}^n$-valued vector of independent random variables. Assume that $||\alpha||=1,$
and that
the function $f$ satisfies
$$ |f(x)-f(x')| \leq d_\alpha(x, x').$$
Let $m$ be a median of $f(X).$ Then for any   $t>0,$
\begin{eqnarray}
\mathbb{P}(f(X) -m  \geq t)     &\leq&  2\,  e^{ -  h(t) } \label{sff1} \\
&\leq& 2\,  e^{ - t^2  }\nonumber
\end{eqnarray}
and
\begin{eqnarray}
\mathbb{P}(m - f(X)\geq t)     &\leq&  2\,  e^{ -  h(t) } \label{sff2} \\
&\leq& 2\,  e^{ - t^2 },\nonumber
\end{eqnarray}
where $h(t)$ is defined by (\ref{sf02})
with $A=\{Y: f(Y)\leq m \}.$
\end{theorem}

Inequalities (\ref{sff1}) and (\ref{sff2}) together  implies the following non-uniform bound between the  medians  and the mean.
\begin{theorem}\label{sf2d} Assume the condition of Theorem \ref{th2}.
Let $m$ be a median of $f(X),$  and let $\mu$ be the mean of $f(X)$. Then
\begin{eqnarray}
| \mu-m |  \leq \big(2  \rho + \sqrt{\pi/2  } \, \big)e^{ -  2 \rho^2  }   < 1.5503,  \label{f15}
\end{eqnarray}
where  $\rho= \mathbb{E}[d_\alpha(X, A)]$ with $A=\{Y: f(Y)\leq m \} $.
\end{theorem}

Inequality  (\ref{f15}) shows that the distance between $\mu$ and $m$ is decaying exponentially  to $0$ as $\rho\rightarrow \infty$.

\section{Extensions to nonnegative functionals }\label{sec3}

Let $ x \in \mathbf{E}^n.$ Denote by  $x^{(i)}=(x_1,...,x_{i-1},x_{i+1},..,x_n) \in \mathbf{E}^{n-1}$ which is obtained by dropping the $i$-th component of $x.$ For each $i\leq n,$ denote $f_i$ a measurable function
 from $\mathbf{E}^{n-1}$ to $\mathbf{R}.$ We have the following extension of
 McDiarmid's inequality (\ref{fsg251}) to a new class of functions. Such class of functions including the self-bounding functions and the
 $(a, b)$-self-bounding introduced by Boucheron, Lugosi and Massart \cite{B00}  and  McDiarmid and Reed \cite{MR0}, respectively.

\begin{theorem}\label{th1ij} Let $X$ be an $\mathbf{E}^n$-valued vector of (not necessarily independent)  random variables. Assume that $||\alpha||=1$, and that for some functions $f_i$,    $i=1,...,n,$ and
all $x \in \mathbf{E}^n,$
\begin{eqnarray}\label{coind01}
0\leq f(x)-f_i(x^{(i)}) \leq \alpha_i.
\end{eqnarray}
Let $\mu$ be the mean  of $f(X).$ Then for any   $t>0,$
\begin{eqnarray}\label{ghjl2}
\mathbb{P}(f(X) -\mu  \geq t)      \leq      e^{ - 2\,t^2 }
\end{eqnarray}
and
\begin{eqnarray}\label{ghjl3}
\mathbb{P}(\mu - f(X)\geq t)      \leq    e^{ - 2\,t^2 } .
\end{eqnarray}
\end{theorem}

Recall that a function $f$ is called $(a, b)$-self-bounding, if
\begin{eqnarray}\label{sfd}
0\leq f(X)-f_i(X^{(i)}) \leq 1
\end{eqnarray}
  and, moreover,  for some $a>0, b\geq0$ and
all $x \in \mathbf{E}^n,$
\begin{eqnarray} \label{conuml}
\sum_{i=1}^{n}\Big( f(x)-f_i(x^{(i)})\Big) \leq af(x)+b.
\end{eqnarray}
In particular,  a $(1, 0)$-self-bounding function is known as a self-bounding function.  It is easy to see that for an $(a, b)$-self-bounding function, condition (\ref{sfd}) implies that $ f(X)/\sqrt{n}$ satisfies  condition (\ref{coind01})   with $\alpha_i=\frac{1}{\sqrt{n}}, i=1,...,n,$ that is
\begin{eqnarray}
0\leq \frac{ f(X)}{\sqrt{n}} - \frac{f_i(X^{(i)})}{\sqrt{n}} \leq  \frac{1}{\sqrt{n}}.
\end{eqnarray}
Thus condition (\ref{sfd}), the inequalities (\ref{ghjl2}) and (\ref{ghjl3}) together implies that for any   $t>0,$
\begin{eqnarray}\label{ghjgl2}
\mathbb{P}(f(X) -\mu  \geq t)      \leq      e^{ - 2\,t^2 /n }
\end{eqnarray}
and
\begin{eqnarray}\label{ghjgl3}
\mathbb{P}( \mu- f(X) \geq t)      \leq      e^{ - 2\,t^2 /n }.
\end{eqnarray}
Notice that the last two inequalities do not assume   condition (\ref{conuml}).

Assume that $X$ is an $\mathbf{E}^n$-valued vector of \emph{independent} random variables.
For any $(a, b)$-self-bounding  function $f,$  McDiarmid and Reed  \cite{MR0} proved that for any $t>0,$
\begin{eqnarray}\label{ffs1}
  \mathbf{P}\Big(\,  f(X)-\mu   \geq t   \Big)
   \leq  \exp\Bigg\{ - \frac{t^2}{2 \big( a \mu  +b   +  a \, t   \big)} \Bigg\}
\end{eqnarray}
and
\begin{eqnarray}\label{ffs2}
  \mathbf{P}\Big(\,  \mu - f(X)  \geq t   \Big)
   \leq  \exp\Bigg\{ - \frac{t^2}{2 \big( a \mu  +b   +    t/3   \big)} \Bigg\}
\end{eqnarray}
 See also  Boucheron, Lugosi and Massart \cite{B00} for the self-bounding functions.
 Theorem  \ref{th1ij}  does not assume condition (\ref{conuml}), and
holds for  \emph{dependent}  random variables. Hence, our inequalities (\ref{ghjl2}) and (\ref{ghjl3}) can be regarded as the  generalizations   of  the  inequalities of Boucheron, Lugosi and Massart \cite{B00} and McDiarmid and Reed \cite{MR0}.

When $X$ is an $\mathbf{E}^n$-valued vector of independent random variables,  Theorem \ref{th1ij}  implies the following   concentration inequalities  around a median.
\begin{theorem}\label{th4}   Let $X$ be an $\mathbf{E}^n$-valued vector of independent random variables. Assume that $||\alpha||=1$, and that  for all $i=1,...,n$ and
all $x \in \mathbf{E}^n,$
\begin{eqnarray}
0\leq f(x)-f_i(x^{(i)}) \leq \alpha_i.
\end{eqnarray}
Let $m$ be a median of $f(X)$,  and  let $\mu = \mathbb{E} [f(X)].$   Then  (\ref{sff1}), (\ref{sff2}) and (\ref{f15}) hold.
\end{theorem}

\section{Proofs of the Theorems} \label{Sec3}
In this section, we devote to the proofs of our theorems.

\noindent\emph{Proof of Theorem \ref{th1}.} Denote by $$f(X)=d_\alpha(X, A) \ \  \ \ \  \textrm{and} \  \ \ \ \  \mu=\mathbb{E}[f(X)].$$
Then it is easy to see that
$$ |f(x)-f(x')| \leq d_\alpha(x, x').$$
By  McDiarmid's  inequality  (\ref{fsg250}), it follows that for any $\lambda\geq 0,$
\begin{eqnarray}\label{sf01f}
\mathbb{E}[e^{ \lambda(f(X) - \mu)}  ]   \leq   \exp\Big\{ \frac{\lambda^2} 8 \Big \}.
\end{eqnarray}
  By  McDiarmid's inequality (\ref{fsg251}),  we get for any $t>0,$
\begin{eqnarray}
\mathbb{P}( \pm ( f(X) - \mu ) \geq t)     \leq   e^{ - 2 \, t^2  }.
\end{eqnarray}
Since $f(X)=0$ if and only if $X \in A,$ we have
\begin{eqnarray}
\mathbb{P}( A)   \leq    \mathbb{P}(   f(X) - \mu    \leq -\mu  )    \leq    e^{ - 2 \, \mu^2  }.
\end{eqnarray}
Thus the last inequality and (\ref{sf01f}) implies that  for any $\lambda\geq 0,$
\begin{eqnarray}
\mathbb{P}( A) \mathbb{E}[e^{ \lambda f(X)}  ]     \leq  \exp \Big\{   \lambda \mu +  \frac{\lambda^2} 8  - 2 \, \mu^2  \Big\}.
\end{eqnarray}
Using Markov's inequality, we get for any $\lambda\geq 0,$
\begin{eqnarray}
\mathbb{P}( A)\mathbb{P}(f(X) \geq t  ) &\leq& \mathbb{P}( A) e^{-\lambda t}\mathbb{E}[e^{ \lambda f(X)}  ]  \nonumber\\
 &\leq & \exp\Big\{ - \lambda t +   \lambda \mu +  \frac{\lambda^2} 8  - 2 \, \mu^2  \Big\}.\label{hknms}
\end{eqnarray}
When $t\geq\mu,$ the right hand side of the last inequality attends its minimum at $\lambda= 4(t- \mu).$ Hence, taking $\lambda= 4(t- \mu)$
in  inequality (\ref{hknms}), we obtain the desired inequality   for  $t\geq\mu$.
When $t\in (0, \mu),$ taking $\lambda= 0$
in  inequality (\ref{hknms}), we obtain the desired inequality  for  $t\in (0, \mu).$
\hfill\qed

\vspace{0.4cm}

\noindent\emph{Proof of Theorem \ref{th2}.} Set $ A=\{ Y \in \mathbf{E}^n: f(Y)\leq m\}. $ Then for each $Y \in A,$
\[
f(X) \leq  f(Y)+ d_\alpha(X, Y).
\]
Minimising over all $Y \in A,$ we have
\[
f(X) \leq  m + d_\alpha(X, A).
\]
Hence,  it is easy to see that $$ \{f(X) -m  \geq t \} \subseteq \{d_\alpha(X, A)  \geq t \}.$$
By Theorem \ref{th1}, for any $t\geq 0,$
\begin{eqnarray*}
\mathbb{P}(f(X) \leq m  )  \mathbb{P}(f(X) -m  \geq t)  &\leq&  \mathbb{P}(A)  \mathbb{P}(d_\alpha(X, A)  \geq t) \\
 &\leq&  e^{ - h(t) }.
\end{eqnarray*}
Since $\mathbb{P}(f(X) \leq m  ) \geq \frac12,$ we obtain the desired inequality (\ref{sff1}). Notice that $-m$ is a media for $-f(X)$ and $\mathbb{P}((-f(X))-(-m)\geq t) =\mathbb{P}(m -f(X)   \geq -t). $  Thus
\begin{eqnarray*}
\mathbb{P}(-f(X) \leq -m  )  \mathbb{P}(m -f(X)   \geq -t)
 \leq   e^{ - h(t)}.
\end{eqnarray*}
By the fact $\mathbb{P}(f(X) \geq m  ) \geq \frac12$ again, we obtain the desired inequality (\ref{sff2}).
 \hfill\qed

\vspace{0.4cm}

\noindent\emph{Proof of Theorem \ref{sf2d}.} It is easy to see that
\begin{eqnarray*}
\mu-m = \mathbb{E}[f(X)-m] \leq \mathbb{E}[(f(X)-m)^+] = \int_0^\infty   \mathbb{P}(f(X) -m  \geq t) dt.
\end{eqnarray*}
By (\ref{sff1}), it follows that
\begin{eqnarray*}
\int_0^\infty   \mathbb{P}(f(X) -m  \geq t) dt &\leq& 2  \int_0^\rho   e^{ -   2 \rho^2}  dt + 2  \int_{\rho}^\infty    e^{ - t^2- (t-2 \rho )^2 }  dt \\
&=&2  \rho   e^{ -   2 \rho^2}   + 2\, e^{ -  2 \rho^2  }   \int_{0}^{\infty}    e^{ -2 \,t^2 }  dt   \\
&=& \Big(2  \rho + \sqrt{\pi/2  } \Big)e^{ -  2 \rho^2  }.
\end{eqnarray*}
Thus $$ \mu-m\leq  \Big(2  \rho + \sqrt{\pi/2  } \Big)e^{ -  2 \rho^2  } . $$
Notice that $-m$ is a media for $-f(X),$ and that $-\mu$ is the means of $-f(X)$. Thus, by (\ref{sff2}),
$$ m-\mu=( -\mu)-(-m)\leq  \Big(2  \rho + \sqrt{\pi/2  } \Big)e^{ -  2 \rho^2  }. $$
In conclusion, it holds
\begin{eqnarray}
| \mu-m | \leq  \Big(2  \rho + \sqrt{\pi/2  } \Big)e^{ -  2 \rho^2  }< 1.5503.
\end{eqnarray}
This completes the proof of Theorem \ref{sf2d}.  \hfill\qed

\vspace{0.4cm}

\noindent\emph{Proof of Theorem \ref{th1ij}.} Denote by $X=(X_1, X_2,...,X_n),$ then $(X_i)_{i=1,..,n}$ is a finite sequence of  random variables.
 Let $(\mathcal{F}_i)_{i=1,..,n}$ be the  natural filtration of the random variables $(X_i)_{i=1,..,n},$
i.e.\ $   \mathcal{F}_i = \sigma\{ X_j,\ 1\leq j \leq i \}. $   Let $f(X)-\mu = \sum_{i=1}^nM_i$ be Doob's martingale decomposition of $f(X),$ where
\begin{eqnarray*}
M_i = \mathbf{E}[f(X) | \mathcal{F}_{i}]-\mathbf{E}[f(X) | \mathcal{F}_{i-1}].
\end{eqnarray*}
By condition (\ref{coind01}), it is easy to see that
\begin{eqnarray*}
M_i &\leq& \mathbf{E}[ \alpha_i + f_i( X^{(i)}) |  \mathcal{F}_{i}]   -\mathbf{E}[  f( X)   | \mathcal{F}_{i-1}]  \nonumber\\
&=&\alpha_i +  \mathbf{E}[  f_i( X^{(i)}) - f(X)     | \mathcal{F}_{i-1}] \label{562s}.
\end{eqnarray*}
Similarly, we have
\begin{eqnarray*}
M_i &\geq& \mathbf{E}[  f_i( X^{(i)}) |  \mathcal{F}_{i}]   -\mathbf{E}[  f( X)   | \mathcal{F}_{i-1}]  \nonumber\\
&=&  \mathbf{E}[   f_i( X^{(i)}) - f( X)     | \mathcal{F}_{i-1}] .\label{572s}
\end{eqnarray*}
Notice that $\mathbf{E}[   f_i( X^{(i)}) - f( X)     | \mathcal{F}_{i-1}]$ is $\mathcal{F}_{i-1}$-measurable. By Azuma's inequality \cite{A67} (see also Corollary 2.7  of \cite{FGL15}), it follows
that  for any $\lambda\geq 0,$
\begin{eqnarray*}
\mathbb{E}[e^{ \lambda(f(X)-\mu)}  ]   \leq \exp\Big\{ \frac{\lambda^2} 8  \sum_{i=1}^n  \alpha_i^2 \Big \} = \exp\Big\{ \frac{\lambda^2} 8 \Big \}.
\end{eqnarray*}
Hence,
\begin{eqnarray}\label{fs}
\mathbb{P}(f(X) -\mu  \geq t)     &\leq&  \inf_{ \lambda \geq 0}  \exp\Big\{-\lambda t+ \frac{\lambda^2} 8 \Big \}= e^{-2 \, t^2  },
\end{eqnarray}
which gives the desired inequality (\ref{ghjl2}).
Notice that  $\mu -f(X) = \sum_{i=1}^n (-M_i)$ and $$-\alpha_i - \mathbf{E}[  f_i( X^{(i)}) - f(X)     | \mathcal{F}_{i-1}]\leq -M_i  \leq- \mathbf{E}[  f_i( X^{(i)}) - f(X)     | \mathcal{F}_{i-1}].$$
By an argument similar to that of (\ref{fs}), we obtained the desired inequality (\ref{ghjl3}).  \hfill\qed

\vspace{0.4cm}

\noindent\emph{Proof of Theorem \ref{th4}.}  Denote by  $x_{(i)}=(x_1,...,x_{i-1}, x_i', x_{i+1},..,x_n) \in \mathbf{E}^{n}$ which is obtained by changing the $i$-th component of $x$ to any given $x_i' \in \mathbf{E}$ such that $x_i' \neq x_i$. Condition (\ref{coind01}) implies that
\begin{eqnarray*}
 f(x) \leq  f_i(x^{(i)}) + \alpha_i \leq f(x_{(i)}) + \alpha_i.
\end{eqnarray*}
By recursion method, the last inequality implies that  for any $x, x' \in \mathbf{E},$
\begin{eqnarray}\label{fdfs01}
 f(x) -   f(x')\leq   d_\alpha(x, x').
\end{eqnarray}
Notice that $ d_\alpha(x, x') =  d_\alpha(x', x).$
By exchanging the places of $x$ and $x',$ we find that (\ref{fdfs01}) is equivalent to that   for any $x, x' \in \mathbf{E},$
\begin{eqnarray}
 |f(x) -   f(x')|\leq   d_\alpha(x, x').
\end{eqnarray}
Thus Theorem \ref{th4} follows by Theorems  \ref{th2} and \ref{sf2d}.\hfill\qed


\end{document}